\documentclass{article} 
\usepackage{amssymb}
\usepackage{amsmath}
\usepackage[all]{xy}  \usepackage[all]{xy}
 \usepackage[psamsfonts]{eucal}
 \usepackage[psamsfonts]{eucal}

\newcommand{\be}{\begin{equation}}
\newcommand{\ee}{\end{equation}}
\newcommand{\bd}{\begin{displaymath}}
\newcommand{\ed}{\end{displaymath}}

\newcommand{\ben}{\begin {enumerate}}
\newcommand{\een}{\end {enumerate}}

\begin {document}

\title{A note on Gorenstein spaces}
\author{Yves F\'elix and Steve Halperin}
\maketitle

\begin{abstract}
Associated with an augmented differential graded algebra $R= R^{\geq 0}$ is a homotopy invariant ${\mathcal T}(R)$. This is a graded vector space, and if $H^0(R)$ is the ground field and $H^{>N}(R)= 0$ then dim$\, {\mathcal T}(R)= 1$ if and only if $H(R)$ is a Poincar\'e duality algebra. In the case of Sullivan extensions $\land W\to \land W\otimes \land Z\to \land Z$ in which dim$\, H(\land Z)<\infty$ we show that 
$${\mathcal T}(\land W\otimes \land Z)= {\mathcal T}(\land W)\otimes {\mathcal T}(\land Z).$$
This is applied to finite dimensional CW complexes $X$ where the fundamental group $G$ acts nilpotently in the cohomology $H(\widetilde{X};\mathbb Q)$ of the universal covering space. If $H(X;\mathbb Q)$ is a Poincar\'e duality algebra and $H(\widetilde{X};\mathbb Q)$ and $H(BG;\mathbb Q)$ are finite dimensional then they are also Poincar\'e duality algebras. \end{abstract}

\vspace{5mm}{\small{\bf Keywords:}
Poincar\'e duality algebras, rational homotopy theory, Gorenstein algebras.}

\vspace{1cm}

In this note we work over an arbitrary ground field $l\!k$. A \emph{Poincar\'e duality algebra} is then a graded algebra $H = \{H^k\}_{0\leq k\leq N}$ such that $H^N = l\!k \,\alpha$ and the pairing
$$\beta \gamma = < \beta, \gamma> \alpha\,, \hspace{1cm} \beta \in H^k, \gamma \in H^{N-k}$$
defines an isomorphism $H^k \stackrel{\cong}{\longrightarrow} \mbox{Hom}(H^{N-k}, l\!k)$, $0\leq k\leq N$. In particular, $H= \mbox{Hom}(\mbox{Hom}(H, l\!k),l\!k)$ is necessarily finite dimensional. 

A \emph{Poincar\'e complex} at $l\!k$ is then a CW complex whose cohomology is a Poincar\'e duality algebra. In this note we develop the properties of path connected spaces satisfying a more general homotopy condition, which reduces to Poincar\'e duality when the cohomology is finite dimensional. As it turns out, this provides additional flexibility which enables applications to fibrations.

This "Gorenstein" condition is defined via the generalization (cf. the Appendix) by Eilenberg and Moore of the classical Ext and Tor functors to the category of modules over a differential graded algebra (dga). More precisely, given a based path connected space $X$ we denote by $C_*(X)$ and $C^*(X)$ the singular chain and cochain complexes for $X$, and by $H_*(X)$ and $H^*(X)$ their respective homology. In particular, the Alexander-Whitney diagonal makes $C^*(X)$ into an augmented dga, and we denote
 $${\mathcal G}(X) = \mbox{Ext}_{C^*(X)}(l\!k, C^*(X)).$$

 \vspace{3mm}\noindent {\bf Definition.}  $X$ is \emph{Gorenstein at $l\!k$} if dim ${\mathcal G}(X)= 1$. In this case ${\mathcal G}(X) = {\mathcal G}(X)^N$ for some $N$, and $N$ is called the \emph{Gorenstein degree} of $X$.
 
 \vspace{3mm} Note that the cap product makes $C_*(X)$ into a right $C^*(X)$-module and that the identification
 $$C^*(X) = \mbox{Hom}_{l\!k}(C_*(X),l\!k)$$
 is an isomorphism of $C^*(X)$-modules. It follows that ${\mathcal G}(X)$ is the dual,
 $${\mathcal G}(X) = \mbox{Hom}_{l\!k} (\mbox{Tor}^{C^*(X)}(l\!k, C_*(X)),l\!k),$$
 and so $X$ is Gorenstein at $l\!k$ if and only if
 $$\mbox{dim Tor}^{C^*(X)}(l\!k, C_*(X)) = 1.$$
 
 \vspace{3mm} The connection with Poincar\'e duality is then provided by
 
 \vspace{3mm}\noindent {\bf Theorem 1.} Suppose an augmented path connected space $X$ satisfies $H_{>N}(X) = 0$, some $N$. Then the following conditions are equivalent:
 \begin{enumerate}
 \item[(i)] $X$ is Gorenstein at $l\!k$.
 \item[(ii)] dim Tor$^{H^*(X)}(l\!k, H_*(X))= 1$.
 \item[(iii)] $H^*(X)$ is a   Poincar\'e duality algebra.
 \end{enumerate}
 
 \vspace{3mm}\noindent {\sl proof:} (i) $\Rightarrow$ (iii) By hypothesis, $C_*(X) = C_*(X)^{\geq -N}$, and so $C_*(X)$ has a minimal $C^*(X)$-semi free resolution $\varphi : P\stackrel{\simeq}{\rightarrow} C_*(X)$. In particular,
 $$\mbox{Tor}^{C^*(X)}(l\!k, C_*(X)) = l\!k \otimes_{C^*(X)}P,$$
 since the differentials in $l\!k \otimes_{C^*(X)}P$ vanish. Since dim Tor$^{C^*(X)}(l\!k, C_*(X))= 1$ it follows that $P$ is the free $C^*(X)$-module on a single generator $a$. Let $z\in C_*(X)$ be the cycle $\varphi (a)$. Then $\varphi (\Phi \cdot a) = \pm z\cap \Phi$ and so the cap product induces an isomorphism of $H^*(X)$-modules
 $$H^*(X) \stackrel{\cong}{\longrightarrow} H_*(X)\,.$$
 Thus   $H^*(X)$    is   a Poincar\'e duality algebra.

 \vspace{2mm} (iii) $\Rightarrow$ (ii). If $H^*(X)$ is a Poincar\'e duality algebra then the cap product makes $H_*(X)$ into a free $H^*(X)$-module with a single generator, the fundamental class. This gives (ii). 
 
 \vspace{2mm} (ii) $\Rightarrow $ (i). this is immediate if Tor$^{C^*(X)}(l\!k, C_*(X))$ is   computed via an Eilenberg-Moore semi free resolution for $C_*(X)$ (cf. the Appendix). 
 
 \hfill$\square$

 The definitions above extend in the obvious way to the category of augmented dga's $R$: we set
 $${\mathcal G}(R) = \mbox{Ext}_R(l\!k, R)\,,$$
 and $R$ is \emph{Gorenstein} if dim$\, {\mathcal G}(R)= 1$. In particular, a quasi-isomorphism $R\stackrel{\simeq}{\rightarrow} S$ of augmented dga's induces an isomorphism ${\mathcal G}(R) \cong {\mathcal G}(S)$ via the isomorphisms
$$\xymatrix{ 
\mbox{Ext}_R(l\!k, R) \ar[rr]^\cong && \mbox{Ext}_R(l\!k, S) && \mbox{Ext}_S(l\!k, S).\ar[ll]^\cong}$$
For simplicity in this category we adopt the notation:
$$-\otimes - = -\otimes_{l\!k} - \,, \hspace{3mm} \mbox{Hom}(-,-) = \mbox{Hom}_{l\!k}(-,-), \hspace{3mm} \mbox{Hom}(M)= \mbox{Hom}_{l\!k}(M,l\!k),$$   and $$\mbox{Hom}_R(M)= \mbox{Hom}_R(M,l\!k),$$
when $R$ is a dga and $M$ is an $R$-module.

This permits the application of dga homotopy theory to the study of Gorenstein spaces in general and Poincar\'e complexes in particular. The relevant definitions and Lemmas from dga homotopy theory are collected in the Appendix at the end of this note.  In particular we introduce the invariants 
$${\mathcal T}(R) = \mbox{Tor}^R(l\!k, \mbox{Hom}(R)) \hspace{3mm}\mbox{and } {\mathcal T}(X) = {\mathcal T}(C^*(X)),$$
for an augmented dga $R$ and a path connected based space $X$. Then we have

 \vspace{3mm}\noindent {\bf Proposition 1.} If $H(R)$ is a graded vector space of finite type, then
 $${\mathcal G}(R) = \mbox{Hom}({\mathcal T}(R)).$$
 In particular, $R$ is Gorenstein if and only if dim$\,{\mathcal T}(R)= 1$
 
 \vspace{3mm}\noindent {\sl proof:}  It follows from the hypotheses that the natural inclusion $R\to \mbox{Hom}(\mbox{Hom}(R))$ is a quasi-isomorphism of $R$-modules, and so
$${\mathcal G}(R) \cong \mbox{Ext}_R(l\!k, \mbox{Hom}(\mbox{Hom}(R)) \cong \mbox{Hom}({\mathcal T}(R)).$$
  \hfill$\square$
 
 \vspace{3mm} In this setting the analogue of Theorem 1 is

\vspace{3mm}\noindent {\bf Theorem 2.}  Suppose $R= R^{\geq 0}$ is a cdga. If $H^0(R)=l\!k$ and $H^{>N}(R)= 0$ then the following conditions are equivalent
\begin{enumerate}
\item[(i)] dim$\, {\mathcal T}(R)= 1$
\item[(ii)] dim$\, {\mathcal T}(H(R))= 1$
\item[(iii)] $H(R) $ is a Poincar\'e duality algebra.
\end{enumerate}
In this case $R$ is Gorenstein.

\vspace{3mm}\noindent {\sl proof.}  The same argument as in the proof of Theorem  1 shows that there is a non-degree preserving isomorphism $\alpha : H(R) \cong H(\mbox{Hom}(R))$. Therefore dim$\, H(R)<\infty$ and $R$ is Gorenstein. Moreover, as in the proof of Theorem 1, as this is an isomorphism of $H(R)$-modules, it follows that $H(R)$ is a Poincar\'e duality algebra. This in turn implies that dim$\, {\mathcal G}(H(R)) = 1$ and so dim$\, {\mathcal T}(H(R))= 1$. Finally using an Eilenberg-Moore semifree resolution for $\mbox{Hom}(R)$ gives that (ii) $\Rightarrow$ (i). \hfill$\square$

\vspace{3mm}\emph{For the rest of this Introduction we restrict to the case $l\!k = \mathbb Q$}. 

\vspace{3mm}
If $X$ is any based, path connected space then $C^*(X)$ is connected by quasi-isomorphisms to the minimal Sullivan model $(\land V,d)$ for $X$ (\cite{FHTI}). This permits the application of Sullivan's theory to establish

\vspace{3mm}\noindent {\bf Theorem 3.} Suppose $F\to X\to Y$ is a fibration of based path connected spaces in which dim$\,H(F)<\infty$,       and $\pi_1(Y)$ acts nilpotently in $H(F)$. Then 
$${\mathcal T}(X) \cong {\mathcal T}(Y)\otimes {\mathcal T}(F).$$
If, additionally,  either $H(Y)$ has finite type, or else
 $H(X)$ has finite type and $H^1(F)= 0$, then  
$${\mathcal G}(X) \cong {\mathcal G}(Y)\otimes {\mathcal G}(F).$$
In this case $X$ is Gorenstein at $\mathbb Q$ if and only if  $H(F)$ is a Poincar\'e duality algebra and $Y$ is Gorenstein at $\mathbb Q$. 

\vspace{3mm}\noindent {\bf Corollary.} Suppose $X$ is a finite dimensional CW complex with fundamental group $G$ and universal covering space $\widetilde{X}$. If $G$ acts nilpotently in each $H^k(\widetilde{X})$ and if $H(X)$ is a Poincar\'e duality algebra,  then each of the conditions  dim$\, {\mathcal T}(\widetilde{X}) = 1$,
  $H(\widetilde{X})$ has finite type,
  $H(BG)$ has finite type, is equivalent to the condition : 
 $H(\widetilde{X})$ is a Poincar\'e duality algebra. When these hold,  $BG$ is Gorenstein at $\mathbb Q$. 

\vspace{3mm}\noindent {\sl proof.} First observe that any of these three conditions implies that dim$\, H(\widetilde{X})<\infty$. In fact, because $\widetilde{X}$ is a finite dimensional complex it follows that for some $N$, $H^{>N}(\widetilde{X})= 0$ and so, if dim$\, {\mathcal T}(\widetilde{X})= 1$, Theorem 2 asserts that dim$\, H(\widetilde{X})<\infty$. Of course if $H(\widetilde{X})$ has finite type then dim$\, H(\widetilde{X})<\infty$. Finally, by hypothesis, dim$\, H(X)<\infty$ and so if $H(BG)$ has finite type, $H(BG)$ and $H(X)$ have Sullivan models of at most countable dimension. Let
$$\land W \to \land W\otimes \land Z\to \land Z$$
be the $\land$-extension determined by the fibration $X\to BG$. Thus $\land W$ and $\land W\otimes \land Z$, are at most countably dimensional, so this also holds for $\land Z$. But by Theorem 5.1 in \cite{FHTII}, $H(\land Z)\cong H(\widetilde{X)}$; Therefore each dim$\, H^p(\widetilde{X})<\infty$, and so dim$\, H(\land Z)= \mbox{dim}\, H(\widetilde{X})<\infty$.  

Now suppose dim$\, H(\widetilde{X})<\infty$. Then Theorem 3 applies and gives 
$${\mathcal T}(X) = {\mathcal T}(BG) \otimes {\mathcal T}(\widetilde{X}).$$
Since $H(X)$ is a Poincar\'e duality algebra, dim$\, {\mathcal T}(X)= 1$ (Theorem 2). Thus dim$\, {\mathcal T}(BG)= \mbox{dim}\, {\mathcal T}(X)= 1$, and so $H(\widetilde{X})$ is a Poincar\'e duality algebra. Moreover, since $Z= Z^{\geq 2}$ and dim$\, H(\land Z)<\infty$ it follows that $Z$ is a graded vector space of finite type. Now Lemma 3.2 in \cite{FHTII} asserts that $H(BG)$ has finite type, and the first  three conditions of the Corollary hold. Since then $H(BG)$ has finite type Proposition 1 implies that $BG$ is Gorenstein at $\mathbb Q$. 

 \hfill$\square$

\vspace{3mm} The hypothesis on $\widetilde{X}$ in the Corollary is essential, and does not always hold even when $BG$ and $X$ are closed manifolds, as the following Example shows.

\vspace{3mm}\noindent {\bf Example.} $X = (S^1\times S^5) \# (S^3\times S^3)$.  In this example, $G = \mathbb Z$ and $BG= S^1$. The universal cover, $\widetilde{X}$ then may be described as follows: Consider first   the connected sum of $(0,1)\times S^5$ with $S^3\times S^3$. The union of this space with $[0,1]\times S^5$ will be denoted by $Z$. We then consider a family of spaces $Z_i$ homeomorphic to $Z$ and indexed by the integers. In each $Z_i$ we denote by $S^5_{i,0}$ the sphere $S^5\times \{0\}$ and by $S^5_{i,1}$ the sphere $S^5\times \{1\}$. Then
$$\widetilde{X} \cong \amalg_i Z_i \,\vert\, (S^5_{i,0}\sim S^5_{i-1,1}).$$
As a CW complex $\widetilde{X}$ has the rational homotopy type of $S^5 \vee (\vee_{i\in \mathbb Z} (S^3\vee S^3))$. The action of $\pi_1(S^1)$ is the identity on $S^5$ and is the translation in the wedge $\vee_{i\in \mathbb Z}\, S^3\vee S^3$.

The model of the projection $X\to S^1$ is
$$(\land u, 0) \to (H^*(X), 0)=  \land (u,z,x,y)/ (uz-xy, ux, uy, zx, zy),$$
with $\vert u\vert = 1$, $\vert z\vert = 5$, $\vert x\vert = \vert y\vert = 3$. A model for the fiber of the model is given by $(H(X)\otimes \land \overline{u}, d)$ with $d\overline{u}= u$. It follows that a basis for the homology of the fiber of the model is given by $1, y, x\otimes \overline{u}^n$ and $y\otimes \overline{u}^n$, $n\geq 1$.  

By a result of Dwyer \cite{Dw} for a fibration $F\to X\to Y$ with nilpotent base, the homology of the fiber of the model is the subalgebra of $H(F)$ in which $\pi_1(Y)$ acts nilpotently. Here $H(F)= H^0(F)\oplus H^5(F)\oplus H^5(F)$, and $H^5(F) = \mathbb Q z$ with trivial action. The space $H^3(F)$ is the space of series $\sum_{n\in \mathbb Z} (a_n t^n \otimes x + b_n t^n \otimes y)$. The invariant elements are the linear combinations $$\alpha \,(\sum t^n\otimes x)+ \beta \,(\sum t^n\otimes y), \hspace{5mm} \alpha, \beta\in \mathbb Q.$$

\vspace{3mm} Theorem 3 is a rational homotopy theory analogue of a theorem of Gottlieb (\cite{Got}) which states that   in a fibration $F\to X\to Y$ of connected spaces which have the homotopy type of finite CW complexes, then $X$ is a Poincar\'e complex if and only if $F$ and $Y$ are. In \cite{gor} the authors introduce the application of dga homotopy theory to establish Theorem 3 under additional finiteness assumptions, and with the hypothesis that the spaces are simply connected. With these additional hypotheses Murillo \cite{mu} extends the result of \cite{gor} to ground fields of arbitrary characteristic.

\section{Fibrations and Sullivan models}

\emph{In this section $l\!k=\mathbb Q$ and all cdga's $R$ satisfy $R= R^{\geq 0}$, and we recall the necessary elements from the theory of Sullivan models.}

\vspace{3mm}
A \emph{Sullivan algebra} is a cdga of the form $\land V$, where the underlying algebra is the free commutative graded algebra on $V = V^{\geq 1}$. The differential is required to satisfy the \emph{Sullivan condition}: $V = \cup_{k\geq 0} V(k)$   where $\{V(k)\}$ is an increasing filtration and $d: V(k+1)\to \land V(k)$. 

More generally, a $\land$-extension is a cdga morphism $R\to R\otimes \land Z$, $x\mapsto x\otimes 1$, in which $Z=Z^{\geq 0}= \cup_{k\geq 0} Z(k)$, and $Z(k)$ is an increasing filtration satisfying $d: Z(k+1) \to R\otimes \land Z(k)$. In particular, $R\otimes \land Z$ is a semifree $R$-module and an augmentation of $R$ extends to the augmentation, $\varepsilon$, of $R\otimes \land Z$ given by $\varepsilon (Z)= 0$. If $Z= Z^{\geq 1}$, this is called a Sullivan extension.

If $H^0(R) = \mathbb  Q$ then $R$ has a $\land$-extension $R\otimes \land U$ with $H(R\otimes \land U)= \mathbb Q$. This is denoted $\overline{R}$ and is called an \emph{acyclic closure} for $R$. In particular, an augmentation $\varepsilon : \overline{R}\to \mathbb Q$ is an $R$-semifree resolution, and thus
\begin{eqnarray}
\label{1.1} 
{\mathcal G}(R) = H(\mbox{Hom}_R(\overline{R},R)), \hspace{5mm}\mbox{and } {\mathcal T}(R) = H(\mbox{Hom}(R)\otimes_R\overline{R}).
\end{eqnarray}

On the other hand, with each path connected space $X$ are associated its Sullivan models; these are Sullivan algebras and are connected by   dga quasi-isomorphisms to $C^*(X)$.  Thus if $\land V$ is a Sullivan model for $X$ then
$${\mathcal G}(\land V) = {\mathcal G}(X) \hspace{3mm}\mbox{and }{\mathcal T}(\land V)= {\mathcal T}(X).$$

Moreover, with a fibration
$$F\to X\to Y$$ of path connected spaces is associated a $\land$-extension
$$\xymatrix{\land W\ar[r]^\lambda & \land W\otimes \land Z \ar[r]^\rho & \land Z}\,, \hspace{1cm} \rho (W)= 0,$$
in which $\land W$ is a Sullivan model for $Y$, $\land W\otimes \land Z$ is a Sullivan model for $X$, and   $Z= Z^{\geq 1}$. 
 
\vspace{3mm}\noindent {\bf Proposition 2.} With  the hypotheses and notation above, suppose in the fibration $F\to X\to Y$ that dim$\, H(F)<\infty$ and $\pi_1(Y)$ acts nilpotently in $H(F)$.   Then
  $\land Z$ is a Sullivan model for $F$, and so
$$ {\mathcal T}(Y)= {\mathcal T}(\land W), \hspace{3mm} {\mathcal T}(X)= {\mathcal T}(\land W\otimes \land Z), \hspace{3mm}\mbox{and }  {\mathcal T}(F)= {\mathcal T}(\land Z).$$
Moreover, if the additional hypotheses of Theorem 3 are satisfied, then each of $H(\land W)$, $H(\land W\otimes \land Z)$, and $H(\land Z)$ are graded vector spaces of finite type. In this case,
$${\mathcal G}(Y) = \mbox{Hom}({\mathcal T}(\land W)), \hspace{3mm} {\mathcal G}(X) = \mbox{Hom}({\mathcal T}(\land W\otimes \land Z)) \hspace{3mm}\mbox{and } {\mathcal G}(F)= \mbox{Hom}({\mathcal T}(\land Z)).$$

\vspace{3mm}\noindent {\sl proof.}   Theorem 5.1 in \cite{FHTII} asserts that $\land Z$ is a Sullivan model for $F$.
If the additional hypotheses are satisfied, either $H(\land W)= H(Y)$ has finite type or else $H(\land W\otimes \land Z)= H(X)$ has finite type and $H^1(F)= 0$. In the first case Proposition 3.8 in \cite{FHTII} asserts that $H(\land W\otimes \land Z)$ has finite type. In the second case, because $H^1(\land Z)= H^1(F)= 0$ and dim$\, H(\land Z)<\infty$,  it follows that $Z$ has finite type. Thus Lemma 3.2 in \cite{FHTII} asserts that $H(\land W)$ has finite type. The final assertion then follows from Proposition 1.\hfill$\square$

\section{Proof of Theorem 3}

Fix a $\land$-extension
$$\xymatrix{R\ar[rr]^\lambda && R\otimes \land Z \ar[rr]^{\rho} && \land Z}$$
in which $R^0= \mathbb Q$ and $Z= Z^{\geq 1}$. Denote $R\otimes \land Z$ by $S$, and let $\overline{R}$, $\overline{\land Z}$ and $\overline{S}$ respectively be  acyclic closures of $R$, $\land Z$ and $S$.

\vspace{3mm}\noindent {\bf Theorem 4.} Suppose in the $\land$-extension above that   dim$\, H(\land Z)<\infty$. Then 
$${\mathcal T}(S) \cong {\mathcal T}(R)\otimes {\mathcal T}(\land Z)$$
In particular, if $H(R)$ has finite type then $R\otimes \land Z$ is Gorenstein if and only if $R$ and $\land Z$ are.

\vspace{3mm} Before embarking on the proof of Theorem 4, consider the special case that    $H(R)$ and $H(\land Z)$ are graded vector spaces of finite type.  
In view of Proposition 2 in the previous section, applied when $R$ is a Sullivan model for $Y$, Theorem 3 is then an immediate consequence of Theorem 4.

\vspace{3mm} For the proof of Theorem 4, we recall that
$${\mathcal T}(R) = H(\mbox{Hom}(R)\otimes_R\overline{R}), \hspace{3mm} {\mathcal T}(S) = H(\mbox{Hom}(S)\otimes_S\overline{S}) \hspace{3mm}\mbox{and } {\mathcal T}(\land Z) = H(\mbox{Hom}(\land Z)\otimes_{\land Z}\overline{\land Z}).$$ 
The proof  is in multiple steps, and we recall from the Appendix that $A^\#$ denotes the underlying graded algebra of a dga $A$. .

\vspace{3mm}\noindent {\bf Proposition 3.} With the hypotheses of Theorem 4,
\begin{enumerate}
\item[(i)] The $S$-modules $\mbox{Hom}_R(S,R)$ is                                                                                                                                                                                       a homotopically semi-free $R$-module.
\item[(ii)] Composition defines a morphism
$$\alpha : \mbox{Hom}_R(S,R)\otimes_R\mbox{Hom}(R)\to \mbox{Hom}(S)$$
of $S$-modules which is also  an $R$-homotopy equivalence. 
\end{enumerate}

\vspace{3mm}\noindent{\sl proof:} Decompose $\land Z$ as the direct sum
\begin{eqnarray}\label{3} 
\land Z = C\oplus \overline{d}C\oplus E\end{eqnarray}
where the differential $\overline{d}$ satisfies $\overline{d}:  C\stackrel{\cong}{\longrightarrow} \overline{d}C$ and $\overline{d}(E)= 0$. Thus
$E\cong H(\land Z)$ and so dim$\, E<\infty$. Then define a morphism
$$\sigma : R\otimes \land Z\to R\otimes \land Z$$
of $R$-modules by setting
$$\sigma (1\otimes x) = 1\otimes x\,, \hspace{1cm} x\in C\oplus E$$
and
$$\sigma (1\otimes \overline{d}x) = d(1\otimes x)\,, \hspace{5mm} x\in \overline{d}C\,.$$

Filtering by the degree in $R$ shows that $\sigma$ is a quasi-isomorphism. It follows that $(R\otimes \land Z)^\#$ is the direct sum of the $R^\#$-modules ($A^\#$ denotes the underlying graded algebra of a dga $A$)
$$(R\otimes \land Z)^\#= (R\otimes C) \oplus \, (R\land d(1\otimes C))\, \oplus (R\otimes E)$$
in which $d: R\otimes C\stackrel{\cong}{\longrightarrow} R\land d(1\otimes C)$. Division by the first two defines a differential $\delta$ in $R\otimes E$ and a surjective quasi-isomorphism
$$\eta : (R\otimes \land Z,d)\to (R\otimes E,\delta)$$
of $R$-modules. Note that  $\rho$ restricts to the retraction $\land Z\to E$ determined by (3). 

Now, because $R\otimes \land Z$ is a $\land$-extension there is an increasing filtration $F_k(\land Z)$, $k\geq 0$, such that $F_0(\land Z)= 0$, $\cup_kF_k(\land Z)= \land Z$, and $d : F_{k+1}(\land Z)\to R\otimes F_k(\land Z)$.  This filtration projects under $\rho$ to a filtration $F_k(E)$ with the corresponding properties with respect to $\delta$. In particular, $(R\otimes E, \delta)$ is $R$-semifree. Since $\eta$ is a quasi-isomorphism of $R$-semifree modules, it is an $R$-homotopy equivalence. This implies in turn that
$$\mbox{Hom}_R(\eta, R) :   \mbox{Hom}_R(R\otimes E,R) \to \mbox{Hom}_R(R\otimes \land Z,R)$$
is also an $R$-homotopy equivalence.

Now since dim$\, E<\infty$, in the filtration $F_k(E)$ some $F_N(E)= E$.   Filter $\mbox{Hom}_R(R\otimes E, R)$ by the submodules $R\otimes G^k$ defined by 
$$G^k = \{ f \in \mbox{Hom}_R(R\otimes E, R)\, \vert \, f(F_{N-k})= 0\}.$$
A simple calculation shows that $G^0= 0$, $G^p\subset G^{p+1}$, $G^N = E$ and
that
$$d : G^{k+1} \to R\otimes G^k\,.$$
It follows that $\mbox{Hom}_R(R\otimes E, R)$ is $R$-semifree and so $\mbox{Hom}_R(R\otimes \land Z, R)$ is a homotopy semifree $R$-module. This establishes (i).

To prove (ii) observe first that $\alpha$ is by definition a morphism of $S$-modules.
Moreover, it follows from (i) that
$$\beta := \mbox{Hom}_R(\eta, R)\otimes_R id : \mbox{Hom}_R(R\otimes E, R)\otimes_R \mbox{Hom}(R) \leftarrow \mbox{Hom}_R(S,R) \otimes_R \mbox{Hom}(R)$$
is an $R$-homotopy equivalence. On the other hand, since dim$\, E<\infty$
$$\mbox{Hom}_R(R\otimes E, R) = \mbox{Hom}(E,R)= \mbox{Hom}(E)\otimes R$$
and so
$$\mbox{Hom}_R(R\otimes E,R)\otimes_R\mbox{Hom}(R)= \mbox{Hom}(E)\otimes \mbox{Hom}(R) = \mbox{Hom}(R\otimes E)\,.$$
Thus we have the commutative diagram
$$\xymatrix{
\mbox{Hom}_R(R\otimes E, R)\otimes_R \mbox{Hom}(R)
 \ar[rr]^\cong  && \mbox{Hom}(R\otimes E) \\
\mbox{Hom}_R(S,R)\otimes_R\mbox{Hom}(R) \ar[u]^\beta\ar[rr]^\alpha && \mbox{Hom}(S)\ar[u]_{\mbox{\scriptsize Hom}(\eta)}.}$$
in which $\beta$ and $\mbox{Hom}(\eta)$ are $R$-homotopy equivalences. If follows that $\alpha$ is an $R$-homotopy equivalence. \hfill$\square$

\vspace{3mm}\noindent{\bf Corollary.} $$\xymatrix{
\overline{S}\otimes_S\mbox{Hom}_R(S,R)\otimes_R\mbox{Hom}(R) \ar[rr]^-{id\otimes_S\alpha} && \overline{S}\otimes_S \mbox{Hom}(S)}$$
is a quasi-isomorphism of $S$-modules.

\vspace{3mm}\noindent {\sl proof:} This follows because $\overline{S}$ is $S$-semi-free and $\alpha$ is a quasi-isomorphism of $S$-modules. \hfill$\square$

\vspace{4mm} Next note that the inclusion $\overline{R}\to \overline{S}$ makes $\overline{S}\otimes_S \mbox{Hom}_R(S,R)$ into an $\overline{R}$-module.

\vspace{3mm}\noindent {\bf Proposition 4.} There is an $\overline{R}$-flat quasi-isomorphism
$$\left[ \, \overline{S}\otimes_S\mbox{Hom}_R(S,R)\right] \otimes_{\overline{R}}  \mathbb Q \otimes \overline{R} \stackrel{\simeq}{\longrightarrow} \overline{S}\otimes_S \mbox{Hom}_R(S,R)\,.$$

\vspace{3mm}\noindent {\sl proof:} Set $T= \overline{R}\otimes_RS$ and observe that 
$$
\overline{S}\otimes_S\mbox{Hom}_R(S,R) = \overline{S}\otimes_TT\otimes_S\mbox{Hom}_R(S,R) = \overline{S}\otimes_T\left( \overline{R}\otimes_R\mbox{Hom}_R(S,R)\right).$$
On the other hand, the inclusion $T\to \overline{S}$ makes $\overline{S}$ into a semifree $T$-module.
 In fact we may write the $T$-module $\overline{S}$ in the form
$$\overline{S}= V\otimes T\,, \hspace{5mm} V = \oplus_{i\geq 0} V(i)\,, \hspace{5mm} V(0)= 0$$
where $d : V(i+1)\to \left( \oplus_{j\leq i}V(j)\right)\otimes T$. Then set
$$\overline{S}(i) = \left[ \oplus_{j\leq i} V(j)\right] \otimes T$$
so that each $\overline{S}(i)$ is a $T$-module and
$$\overline{S}(i+1)/\overline{S}(i) = (V(i+1)\otimes T, id\otimes d)\,.$$

Then we may write
$$
\overline{S}\otimes_S\mbox{Hom}_R(S,R)  = \cup_i \overline{S}(i) \otimes_T (\overline{R}\otimes_R \mbox{Hom}_R(S,R))\,. 
$$
Since $\mbox{Hom}_R(S,R)$ is a homotopically semi-free $R$-module and $\overline{R}$ is $R$-semi-free it follows that $\overline{R}\otimes_R \mbox{Hom}_R(S,R)$ is homotopically $\overline{R}$-semi-free. In particular, it is $\overline{R}$-flat. For simplicity write
$$M:= \overline{R}\otimes_R\mbox{Hom}_R(S,R)\,.$$

Now  observe that each
$$\overline{S}(i+1)/\overline{S}(i)\otimes_TM = \left( V(i+1)\otimes T\right)\otimes_TM= V(i+1)\otimes M$$
is an $\overline{R}$-flat module. It follows by induction on $i$ that each $\overline{S}(i)\otimes_TM$ is $\overline{R}$-flat and hence $\overline{S}\otimes_TM$ is $\overline{R}$-flat.
Moreover, since $\varepsilon : \overline{R}\to \mathbb Q$ is a quasi-isomorphism of $\overline{R}$-modules it then follows that
$$\varepsilon (i) : \overline{S}(i)\otimes_TM \to (\overline{S}(i)\otimes_TM)\otimes_{\overline{R}}\mathbb Q \hspace{5mm}\mbox{and } \varepsilon_M : \overline{S}\otimes_TM \to \overline{S}\otimes_TM\otimes_{\overline{R}} \mathbb Q$$
are surjective quasi-isomorphisms.

The next step  is to construct maps of complexes
$$\sigma (i) : (\overline{S}(i)\otimes_TM)\otimes_{\overline{R}}\mathbb Q \to \overline{S}(i)\otimes_TM\,, \hspace{5mm} \sigma : (\overline{S}\otimes_TM) \otimes_{\overline{R}}\mathbb Q \to \overline{S}\otimes_TM$$
such that $\varepsilon (i) \circ \sigma (i) = id$, $\sigma (i+1)$ extends $\sigma (i)$ and $\sigma = \varinjlim_i \sigma (i)$. Suppose $\sigma (i)$ is constructed. Since $\varepsilon (i+1)$ is a surjective quasi-isomorphism it has a right inverse $\omega$. Then
$$\sigma (i)-\omega : \overline{S}(i)\otimes_TM\otimes_{\overline{R}} \mathbb Q \to \mbox{ker}\, \varepsilon (i+1).$$
Because $H(\mbox{ker}\, \varepsilon (i+1))= 0$, there is a map $\gamma : \overline{S}(i)\otimes_TM\otimes_{\overline{R}}\mathbb Q \to \mbox{ker}\, \varepsilon (i+1)$ such that 
$$\sigma (i)-\omega = d\circ \gamma + \gamma \circ d.$$
Extend $\gamma$ to a map
$$\gamma : \overline{S}(i+1)\otimes_TM\otimes_{\overline{R}}\mathbb Q \to \mbox{ker}\, \varepsilon (i+1)$$
and set $\sigma (i+1) = \omega + d\circ \gamma + \gamma \circ d$. Then set $\sigma = \varinjlim_i \sigma (i)$. Since each $\varepsilon (i)$ is a quasi-isomorphism so are each $\sigma (i)$ and $\sigma$.

Further, $\sigma (i)$ and $\sigma$ extend uniquely to morphisms of $\overline{R}$-modules
$$\tau (i) : \overline{S}(i)\otimes_TM\otimes_{\overline{R}}\mathbb Q\otimes  \overline{R} \to \overline{S}(i)\otimes_TM$$
and
$$\tau : \overline{S}\otimes_TM\otimes_{\overline{R}}\mathbb Q\otimes \overline{R}\to \overline{S}\otimes_TM.$$
 The commutative diagram
$$
\xymatrix{
\overline{S}(i)\otimes_TM\otimes_{\overline{R}}\mathbb Q \otimes \overline{R} \ar[rr]^{\tau (i)} && \overline{S}(i)\otimes_TM\\
\overline{S}(i)\otimes_TM\otimes_{\overline{R}}\mathbb Q \ar[u]^{\lambda (i)}_\simeq \ar[rru]^{\sigma (i)}
}\hspace{1cm} \lambda (i) \Phi = \Phi \otimes 1$$
shows that each $\tau (i)$ is a quasi-isomorphism. Thus so is $\tau$.

We now show by induction that each $\tau (i)$ is $\overline{R}$-flat. In fact, $\tau (i)$ and $\tau (i+1)$ induce a quotient quasi-isomorphism
$$\nu : \overline{S}(i+1)/\overline{S}(i) \otimes_TM\otimes_{\overline{R}}\mathbb Q\otimes \overline{R} \to \overline{S}(i+1)/\overline{S}(i)\otimes_TM,$$
which may be identified as a quasi-isomorphism
$$V(i+1)\otimes M\otimes_{\overline{R}}\mathbb Q \otimes \overline{R} \stackrel{\simeq}{\longrightarrow} V(i+1)\otimes M$$
from a semi-free $\overline{R}$-module to a homotopically semi-free $\overline{R}$-module. Such a morphism is a homotopy equivalence and hence is $\overline{R}$-flat. Thus if $\tau (i)$ is $\overline{R}$-flat so is $\tau (i+1)$, and thus so is $\tau = \varinjlim_i \tau (i)$.

But
$$\overline{S}\otimes_TM= \overline{S}\otimes_T(\overline{R}\otimes_R \mbox{Hom}_R(S,R)) = \overline{S}\otimes_S \mbox{Hom}_R(S,R)$$
and so the Proposition is proved. 

\hfill$\square$

 Recall that if $B\subset A$ is a sub dga and $B^0= \mathbb Q$ then we write
 $$A//B := A\otimes_B\mathbb Q\,.$$
 In particular, if $M$ is any $S$-module then the projection
 $$\rho_M : M\to M\otimes_R\mathbb Q$$
 induces an $S//R$-module structure in $M\otimes_R\mathbb Q$. Similarly, if $N$ is a $T$-module then $\mathbb Q\otimes_{\overline{R}}N$ inherits a $T//\overline{R}$-module structure. But
 $$T//\overline{R} = (S\otimes_R\overline{R})\otimes_{\overline{R}}\mathbb Q = S//R$$
 and so $\mathbb Q\otimes_{\overline{R}}N$ is an $S//R$-module.
 
 \vspace{3mm}\noindent {\bf Proposition 5.} 
 If $N$ is a semi-free $T$-module and $M$ is an $S$-module then
 $$\xymatrix{
(N\otimes_SM)\otimes_{\overline{R}}\mathbb Q =  (\mathbb Q\otimes_{\overline{R}}N)\otimes_SM \ar[rrr]^{id\otimes \rho_M} &&&  (\mathbb Q\otimes_{\overline{R}}N)\otimes_{S//R} (M \otimes_R\mathbb Q)}$$
 is an isomorphism.

 \vspace{3mm}\noindent {\sl proof:}  An exact sequence and direct limit argument reduces the Proposition to the case $N=T$. In this case $id\otimes \rho_M$ has the form
 $$S//R \otimes_SM \to S//R\otimes_{S//R}(M\otimes_R\mathbb Q) = M\otimes_R\mathbb Q\,.$$
 But
 $$S//R \otimes_SM = \mathbb Q \otimes_RS\otimes_SM = \mathbb Q \otimes_RM,$$
 which identifies $id\otimes \rho_M$ as the identity. \hfill$\square$

 \vspace{3mm}\noindent {\sl proof of Theorem 4:} As above we denote $S:= R\otimes \land Z$ and $T = \overline{R}\otimes_RS$.

 Since $\overline{S}$ is $S$-semi-free, Proposition 3 provides a quasi-isomorphism
 $$\overline{S}\otimes_S\mbox{Hom}_R(S,R)\otimes_R\mbox{Hom} (R)\simeq \overline{S}\otimes_S \mbox{Hom} (S)\,.$$
 Write the left hand side as
 $$\overline{S}\otimes_S \mbox{Hom}_R(S,R) \otimes_{\overline{R}} (\overline{R}\otimes_R\mbox{Hom}(R))\,.$$
 Proposition 4 then provides an $\overline{R}$-flat quasi-isomorphism
 $$\left[ \overline{S}\otimes_S\mbox{Hom}_R(S,R)\right] \otimes_{\overline{R}} \mathbb Q\otimes \overline{R}\stackrel{\simeq}{\longrightarrow} \overline{S}\otimes_S \mbox{Hom}_R(S,R)\,.$$
 From this we obtain a quasi-isomorphism
 $$ \left[\overline{S}\otimes_S\mbox{Hom}_R(S,R)\right] \otimes_{\overline{R}} (\overline{R}\otimes_R \mbox{Hom}(R)) \simeq \left[ \, [ \overline{S}\otimes_S\mbox{Hom}_R(S,R)\otimes_{\overline{R}}\mathbb Q\, \right] \otimes (\overline{R}\otimes_R\mbox{Hom}(R))\,.$$
 
 It remains to show that
 $$\left[ \overline{S}\otimes_S \mbox{Hom}_R(S,R)\right]\otimes_{\overline{R}}\mathbb Q \simeq \overline{\land Z}\otimes_{\land Z} \mbox{Hom}(\land Z)\,.$$
Since $\overline{S}$ is semi-free, Proposition 5 provides a quasi-isomorphism
$$ [\overline{S}\otimes_S\mbox{Hom}_R(S,R)]\otimes_{\overline{R}}\mathbb Q \simeq \overline{S}//\overline{R} \otimes_{S//R} (\mbox{Hom}_R(S,R)\otimes_R\mathbb Q)\,.$$

Moreover, there is a natural inclusion of $S$-modules
$$\mbox{Hom}_R(S,R)\otimes R \to \mbox{Hom}_R(S,R)$$
which we show is an $R$-homotopy equivalence. In fact, in the proof of Proposition 3 we constructed a  homotopy equivalence of $R$-modules,
$$\eta : S\to R\otimes E$$
 in which dim$\, E<\infty$. This yields the commutative diagram
$$\xymatrix{
\mbox{Hom}_R(S, \mathbb Q)\otimes R \ar[rr] && \mbox{Hom}_R(S,R)\\
\mbox{Hom}_R(R\otimes E, \mathbb Q) \otimes R \ar[u]\ar[rr] && \mbox{Hom}_R(R\otimes E, R)\ar[u]}$$
in which the vertical arrows are homotopy equivalences of $R$-modules. But the bottom horizontal arrow is an isomorphism because dim$\, E<\infty$. it follows that
$$\mbox{Hom}_R(S,\mathbb Q)\otimes R \to \mbox{Hom}_R(S,R)$$
is an $R$-homotopy equivalence. Applying $\otimes_R\mathbb Q$ therefore gives a quasi-isomorphism
$$\mbox{Hom}_R(S, \mathbb Q) \stackrel{\simeq}{\longrightarrow} \mbox{Hom}_R(S,R)\otimes_R\mathbb Q$$
of $S$-modules.

Now, $\overline{S}/\!/\overline{R}$ is $E/\!/R$-semi-free, and it follows that
$$\overline{S}/\!/\overline{R}\otimes_{S/\!/R} (\mbox{Hom}_R(S, \mathbb Q))\stackrel{\simeq}{\longrightarrow} \overline{S}//\overline{R} \otimes_{S/\!/R} (\mbox{Hom}_R(S,R)\otimes_R\mathbb Q)\,.$$
But $\mbox{Hom}_R(S, \mathbb Q)= \mbox{Hom}(S/\!/R)$. Thus these quasi-isomorphisms combine to give
$$\overline{S}\otimes_S\mbox{Hom}_R(S, R)\otimes_{\overline{R}}\mathbb Q \simeq \overline{S}/\!/\overline{R} \otimes_{S//R} \mbox{Hom}(S//R)\,.$$
But $\overline{S}/\!/\overline{R}= \overline{\land Z}$ and $S/\!/R= \land Z$, which establishes the isomorphism of the Theorem.

Finally, as observed at the start of this section, if $H(R)$ has finite type so does $H(S)$. In this case ${\mathcal G}(R)$, ${\mathcal G}(S)$ and ${\mathcal G}(\land Z)$ are respectively the duals of ${\mathcal T}(R)$, ${\mathcal T}(S)$ and ${\mathcal T}(\land Z)$ This gives the last assertion of the Theorem.   \hfill$\square$.

 \section{Sullivan extensions} \emph{In this section $l\!k = \mathbb Q$}. 
 
 As recalled in Section 1 a fibration $F\to X\to Y$ of path connected spaces determines a Sullivan extension
 $$\land W\to \land W\otimes\land Z\to \land Z$$
 in which $\land W$ is a Sullivan model for $F$ and $\land W\otimes \land Z$ is a Sullivan model for $X$. Here $\land W$ and $\land Z$ may be chosen to be minimal models: $d\colon W\to \land^{\geq 2}W$ and $\overline{d}: Z\to \land^{\geq 2}Z$. Moreover, if $\pi_1(Y)$ acts nilpotently in each $H^p(F)$ and each dim $H^p(F)<\infty$ then $\land Z$ is a Sullivan model for $F$. 
 
 Now consider the special case of the fibration $\widetilde{X}\to X\to BG$ described in the Corollary to Theorem 3, in which $X$ is a finite dimensional CW complex. Assume further that    $H(X)$ and $H(\widetilde{X})$ are Poincar\'e duality algebras, and   dim$\, H(BG)<\infty$. 
 
 \vspace{3mm}\noindent {\bf Theorem 5.} Suppose in the minimal Sullivan model $\land W$ of $BG$ that $W^{\geq 2}$ has finite type. Then, with the hypotheses above, $H(\land W)$ and $H(\land W^{\geq 2})$ are Poincar\'e duality algebras. 
 
 \vspace{3mm}\noindent {\bf Remark.} If $\land V$ is the minimal Sullivan model of a path connected space $X$ then the sub dga $\land V^1$ is the Sullivan analogue of $B(\pi_1(X))$, and indeed the minimal Sullivan model of $B(\pi_1(X))$ has the form $\land V^1\otimes \land W^{\geq 2}$. But it may well happen that $W^{\geq 2}\neq 0$.
 
 \vspace{3mm}\noindent {\sl proof ot Theorem 5.} The acyclic closure of $\land W^1$ has the form $\land W^1\otimes \land U$ with $U$ concentrated in degree $0$. Since dim$\, H(\land W)<\infty$ it follows that
 $$H(\land W^{\geq 2}) = H(\land W\otimes_{\land W^1}\land W\otimes \land U) = H(\land W\otimes \land U)$$
 satisfies $H^{>N}(\land W^{\geq 2})= 0$ for some $N$. Since $W^{\geq 2}$ is assumed to have finite type this implies that dim$\, H(\land W^{\geq 2})<\infty$.
 
 Now by Theorem 4,
 $${\mathcal T}(\land W) = {\mathcal T}(\land W^1) \otimes {\mathcal T}(\land W^{\geq 2}).$$
 But since $H(\land W)$ is is a Poincar\'e duality algebra,  this implies that ${\mathcal T}(\land W^1)$ and ${\mathcal T}(\land W^{\geq 2})$ are 1-dimensional. In particular, $H(\land W^{\geq 2})$ is a Poincar\'e duality algebra. Finally, since dim$\, H(\land W)<\infty$, Proposition 9.6 in \cite{FHTII} asserts that $H^{>M}(\land W^1)= 0$ for some $M$. Now Theorem 2 asserts that $H(\land W')$ is a Poincar\'e duality algebra. \hfill$\square$
 
 \section{An extension of Theorem 3}
 
 \emph{Here again, $l\!k = \mathbb Q$} and we consider a fibration
 $$F\to X\to Y$$ of path connected spaces with corresponding Sullivan extension
 $$\land W\to \land W\otimes \land Z\to \land Z.$$
Then  $H(\land W)\cong H(Y)$ and $H(\land W\otimes \land Z) \cong H(X)$, but if $\pi_1(Y)$ does not act nilpotently in $H(F)$ it may happen that 
 $H(\land Z)$ and $H(F)$ are not isomorphic.
 
 In this setting we denote by $U\subset H(F)$ the subalgebra of elements on which $\pi_1(Y)$ acts nilpotently, and we recall that for any $\pi_1(Y)$-module $M$, $H(Y;M)$ denotes the cohomology of $Y$ with local coefficients in $M$.
 
 \vspace{3mm}\noindent {\bf Theorem 6.} Suppose with the notation above that   $H(B)$ has finite type and that dim$\, U<\infty$. If 
 $$H(Y;U)\stackrel{\cong}{\longrightarrow} H(Y,H(F)),$$
 then $${\mathcal T}(X) \cong {\mathcal T}(Y) \otimes {\mathcal T}(U).$$
 In particular, if $X$ is Gorenstein at $\mathbb Q$ then $U$ is a Poincar\'e duality algebra.
 
 \vspace{3mm}\noindent {\sl proof:} Theorems 1 and 3 in \cite{FT} imply that $U\cong \land Z$, and so Theorem 6 follows from Theorem 4. \hfill$\square$

 \section{Appendix: Differential algebra}

We establish the following conventions   in the category of modules over   differential graded algebras (dga) $A$ defined over an arbitrary ground field $l\!k$. 
  The differentials are usually denoted by $d$, but are normally suppressed from the notation: if $M$ is a differential graded object then $M^{\sharp}$ denotes the underlying graded object.  
 If $A$ is a commutative dga (cdga) then  
left $A$-modules determine the corresponding right $A$-modules via
$$
x.a=(-1)^{deg\,a\cdot deg\, x}a.x,   a\in A. 
$$
These are identified when $A$ is a commutative dga (cdga). 

 The functors $Hom_{l\!k}(-,l\!k)$ and $-\otimes_{l\!k}-$ are simply denoted $Hom(-)$ and $-\otimes-$.
  If $A^0= l\!k$ and $M$ is an $A$-module then $M\otimes_A l\!k$ is sometimes denoted by $M//A$.
  
 If $\varphi:R\rightarrow A$ is a  dga morphism and $M$ is an R-module then $Hom_{R}(A,M)$ is assigned the A- and R- module structures given by
$$
f.a(x)=f(a.x) \hbox{\ \ and\ \ } f.b(x)=(-1)^{deg\,b\cdot deg\,x}f\,(x).b, \hspace{5mm} x\in M, a\in A, b\in R.
$$
Thus this R-module structure coincides with that obtained from the A-module structure via precomposition with $\varphi$.

 Two A-modules $M$ and $N$ have the same homotopy type if there are morphisms $\varphi:M\rightarrow N$ and $\psi:M\leftarrow N$ such that both composites are homtopic to the respective identities via $A^{\sharp}$-linear maps.

\vspace{3mm}
Next recall that an A-module P is \textit{A-semifree} if it can be equipped with an increasing filtration by A-modules $F_{k}, k\geq 0,$ such that $F_{0}=0, P=\bigcup_{k}F_{k}$ and there are A-module isomorphisms 
$$
F_{k+1}/F_{k}\cong (A\otimes S_{k}, d\otimes id), k\geq 0.$$
  In this case we have an isomorphism $A\otimes V\stackrel{\cong}{\rightarrow} P$ in which $V = \oplus_{k\geq 0} S_k$ and $d : S_{k+1}\to A\otimes (\oplus_{j\leq k} S_j)$.

A \emph{semifree resolution} of an $A$-module $M$ is a quasi-isomorphism $P\stackrel{\simeq}{\rightarrow} M$ from a semi free $A$-module. These always exist and are unique up to quasi-isomorphism (\cite[Chap.6]{FHTII}). Given a semi free resolution $P\stackrel{\simeq}{\to} M$ of an $A$-module,  filtering by the $F_k$ determines a spectral sequence whose $E_1$-term has the form
$$\to H(A)\otimes S_{k+1} \to H(A)\otimes S_k \to \dots \to H(A)\,.$$
The semi free resolution is called an \emph{Eilenberg-Moore resolution} if this sequence is a resolution of the $H(A)$-module $H(M)$. Every $A$-module has an Eilenberg-Moore resolution (\cite[Prop 20.11]{FHTI}).

The Eilenberg-Moore generalizations of Ext and Tor to dga-modules are defined by
$$\mbox{Ext}_A(M,-) = H(\mbox{Hom}(P,-))\hspace{1cm}\mbox{and } \mbox{Tor}Â(M,-)= H(P\otimes_A-),$$
where $P$ is any semi free resolution of $M$. 
In particular, if $N$ is a second $A$-module then
$$ 
\mbox{}\hspace{3cm} \mbox{Ext}_A(M, \mbox{Hom}(N)) = \mbox{Hom}(\mbox{Tor}_A(M,N)). \hspace{4cm} \mbox{(A1)}
$$

Now consider the special case of a dga $A$ in which $A= A^{\geq 0}$ and $H^0(A)=l\!k$. Let $S\subset A^1$ be a direct summand of $d(A^0)$ and set 
$$B^0= l\!k\,, \hspace{5mm} B^1= S,\hspace{5mm}\mbox{and } B^k= A^k\,, \hspace{3mm} k\geq 2.$$
This defines a sub dga $B\subset A$ for which the inclusion is a quasi-isomorphism. Since $A\otimes_B $ converts semi-free $B$-modules to semi free $B$-modules the standard construction for $B$-modules gives

\vspace{3mm}\noindent {\bf Lemma A.1.} {\sl Suppose a dga $A$ satisfies $A= A^{\geq 0}$ and $H^0(A)= l\!k$. If an  $A$-module $M$ satisfies $H(M)= H^{\geq k}(M)$, some $k\in \mathbb Z$, then $M$ has an $A$-semifree resolution $A\otimes V\stackrel{\simeq}{\rightarrow} M$ such that $d(V)\subset A^{\geq 1}\otimes V$.}

\vspace{3mm}\noindent {\bf Definition.} A semi free $A$-module $A\otimes V$ is called minimal if $d(V) \in A^{\geq 1}\otimes V$. 

\vspace{3mm}\noindent {\bf Lemma A.2.}  {\sl Suppose $\varphi : M\to N$ is a morphism of $A$-modules in which $A$ is augmented to $l\!k$, $A= A^{\geq 0}$, and $H^0(A)=l\!k$. If $M$ and $N$ admit minimal semi free resolutions, then the following conditions are equivalent:
\begin{enumerate}
\item[(i)] $\varphi$ is a quasi-isomorphism
\item[(ii)]$ \mbox{Tor}_A(l\!k, \varphi)$ is an isomorphism.
\end{enumerate}}

\vspace{3mm}\noindent {\sl proof:} That (i) $\Rightarrow$ (ii) is standard (see eg. \cite[Chap.6]{FHTII}). To see that (ii) $\Rightarrow$ (i) replace $A$ by the dga $B$ described above to reduce to the case $A^0= l\!k$. Then lift $\varphi$ up to homotopy to a morphism $\psi : P\to Q$ between minimal semi free resolutions for $M$ and $N$. Since the differentials in $l\!k\otimes_AP$ and $l\!k\otimes_AQ$ vanish, $\psi$ induces an isomorphism $l\!k\otimes_AP\stackrel{\simeq}{\longrightarrow} l\!k \otimes_AQ$.

Now write $P= A\otimes V$ and $Q= A\otimes W$, and filter by the subspaces $A^{\geq k}\otimes -$. Then the associated graded map induced by $\psi$ is an isomorphism, and so $\psi$ itself is an isomorphism and $\varphi$ is a quasi-isomorphism. \hfill$\square$

\vspace{3mm}
More generally, we say an $A$-module $M$ is \emph{homotopy semifree} 
  if it has the homotopy type of a semifree A-module.

\vspace{3mm}\noindent {\bf Lemma A.3} \begin{enumerate}
\item[(i)] Suppose $M$ is a homotopy semifree A-module. Then $M\otimes_{A}-$ preserves quasi-isomorphisms.
 \item[(ii)] A quasi-isomorphisms between semi-free $A$-modules   is a homotopy equivalence. 
 \item[(iii)]  If $P$ is any $A$-module and if the A-modules M and N have the same homotopy type then so do the A-modules $Hom_{A}(M,P)$ and $Hom_{A}(N,P)$. 
 \item[(iv)] A quasi-isomorphism $\varphi : M\to N$ from a semi-free $A$-module to a homotopy semifree $A$-module is a homotopy equivalence.\end{enumerate}

\vspace{3mm}\noindent {\sl proof:} Immediate from the definitions. \hfill$\square$

\vspace{3mm}Next, suppose $A$ is a cdga and $A^0=\mathbb Q$. Then (\cite[Chap.3]{FHTII}) there  is a cdga morphism
$$\iota_A : A\to \overline{A}$$
such that $\iota_A$ is injective, $\overline{A}$ is $A$-semifree, and $H(\overline{A})=\mathbb Q$. The cdga $\overline{A}$ is called an \emph{acyclic closure} for $A$.

If $A$ is a cdga then an $A$-module $M$ is \emph{flat} if $M\otimes_A-$ preserves quasi-isomorphisms. Also
   a quasi-isomorphism $\varphi : M\stackrel{\simeq}{\to} N$ of $A$-modules is \emph{flat} if for any $A$-module, $P$, $\varphi\otimes_Aid_P$ is a quasi-isomorphism. 
   
\vspace{3mm}\noindent {\bf Lemma A.4} \mbox{}\\
\begin{enumerate}
\item[(i)] If $0\to M(1)\to M(2)\to M(3)\to 0$ is an exact sequence of $A$-modules and $M(1)$ and $M(3)$ are flat, then so is $M(2)$.
\item[(ii)] The direct limit of flat $A$-modules is flat.
\item[(iii)] A homotopically semifree $A$-module  is  flat.
\end{enumerate} 
 
\vspace{2mm}\noindent {\sl proof}: (i) and (ii) are immediate. By \cite{FHTI} semifree $A$-modules are flat, and it follows that so are homotopy semifree $A$-modules. \hfill$\square$

\vspace{3mm}\noindent {\bf Lemma A.5} \mbox{}\\
\begin{enumerate}
\item[(i)] The direct limit of flat quasi-isomorphisms of $A$-modules is a flat quasi-isomorphism.
\item[(ii)] Suppose
$$\xymatrix{
0 \ar[r] &  M(1)\ar[r]\ar[d]^{\varphi (1)} & M(2)\ar[r]\ar[d]^{\varphi (2)} & M(3) \ar[r]\ar[d]^{\varphi (3)} & 0\\
0 \ar[r] & N(1)\ar[r] & N(2)\ar[r] & N(3)\ar[r] & 0}$$
is a row exact commutative diagram of morphisms of $A$-modules. If $\varphi (1)$ and $\varphi (3)$ are $A$-flat so is $\varphi (2)$.\end{enumerate}

   \vspace{1cm} Institut de Recherche en Math\'ematique et en Physique, Universit\'e Catholique de Louvain, 2, Chemin du Cyclotron, 1348 Louvain-La-Neuve, Belgium, yves.felix@uclouvain.be
   
   \vspace{3mm} Department of Mathematics, Mathematics Building, University of Maryland, College Park, MD 20742, United States, shalper@umd.edu

 \end{document}